\documentclass[12pt,oneside,reqno]{amsart}
\usepackage{graphicx}

\graphicspath{ {figs/} }

\usepackage{pict2e}
\usepackage{amssymb}
\usepackage{amsthm}                
\usepackage[margin=0.9in]{geometry}  

\usepackage{epstopdf}
\usepackage{amsmath}
\usepackage{graphicx}
\usepackage{subfigure}
\usepackage{latexsym}
\usepackage{amsmath}
\usepackage{amsfonts}
\usepackage{amssymb}
\usepackage{tikz}
\usepackage{mathdots}
\usepackage{comment}
\usepackage{float}
\usepackage[mathscr]{euscript}
\usepackage{enumerate}
\usepackage{epsfig}
\usepackage { hyperref }
\usepackage { graphics }
\usepackage { graphicx }
\usepackage{arydshln}
\usepackage{cite}

\usetikzlibrary{arrows.meta}
\usetikzlibrary{knots}
\usetikzlibrary{hobby}
\usetikzlibrary{arrows,decorations.markings}

\usepackage[utf8]{inputenc}
\usepackage{longtable}
\usepackage[lite]{amsrefs}

\renewcommand{\PrintDOI}[1]{\href{http://dx.doi.org/\detokenize{#1}}{doi: \detokenize{#1}}%
	\IfEmptyBibField{pages}{, (to appear in print)}{}}

\theoremstyle{definition}
\newtheorem{theorem}{Theorem}[section]

\newtheorem{corollary}[theorem]{Corollary}
\newtheorem{proposition}[theorem]{Proposition}

\theoremstyle{definition}
\newtheorem{definition}[theorem]{Definition}
\newtheorem{example}[theorem]{Example}
\theoremstyle{remark}

\numberwithin{equation}{section}

\numberwithin{equation}{section}

\def\utr{\, \underline{\triangleright}\, }
\def\otr{\, \overline{\triangleright}\, }
\def\ud{\, \underline{\bullet}\, }
\def\od{\, \overline{\bullet}\, }

\title{Singquandles, Psyquandles and Singular Knots: A Survey}

\author{Jose Ceniceros}
\address{Hamilton College, Clinton, NY }
\email{jcenicer@hamilton.edu}

\author{Indu R. Churchill}
\address{State University of New York at Oswego, Oswego, NY }
\email{indurasika.churchill@oswego.edu}

\author{Mohamed Elhamdadi}
\address{University of South Florida, Tampa, FL }
\email{emohamed@math.usf.edu}

\author{Mustafa Hajij}
\address{Santa Clara University,Santa Clara, CA}
\email{hajij@scu.edu}

\date{}
\subjclass[2020]{Primary 57K12, 05C38; Secondary 05A15}
\keywords{Quandles, singular Knots and Links}
\dedicatory{}

\begin{document}
\maketitle 
	
\begin{abstract}

 In this short survey we review recent results dealing with algebraic structures (quandles, psyquandles, and singquandles) related to singular knot theory. We first explore the singquandles counting invariant and then consider several recent enhancements to this invariant. These enhancements include a singquandle cocycle invariant and several polynomial invariants of singular knots obtained from the singquandle structure. We then explore psyquandles which can be thought of as generalizations of oriented signquandles, and review recent developments regarding invariants of singular knots obtained from psyquandles. 
\end{abstract}

\tableofcontents

\section{Introduction}
The goal of this article is to summarize some recent works on some algebraic structures related to singular knot theory. 
In this article ``A quick trip through knot theory" \cite{Fox}, Fox introduced in 1961 a diagrammatic definition of colorability of a knot by the set of integers modulo $n$.  This notion of colorability is obviously one of the simplest invariant of knots, see for example \cite{Przy}.  For a natural number $n>1$, a diagram of a knot is $n$-colorable if at every crossing, the sum of the colors of the under-arcs is twice the color of the over arc modulo $n$.  This integer valued invariant has been generalized to give quandle coloring-counting invariant and applied to many knotted objects such as tangles, braids, knots etc, see \cite{EN}.

In early 80s,  the notion of a quandle was introduced independently by Joyce \cite{Joyce} and Matveev \cite{Matveev} and used it to construct representations of the braid groups and thus invariants of knots and links.  Since then there has been lot of work on quandles and their ramifications by both topologists and algebraists.

The theory of singular knots was introduced in the 90s by V. A. Vassiliev in \cite{V} as an extension of classical knot theory.  Since then several classical knot invariants have been extended to singular knots.  Birman in \cite{Birman} investigated the braid theory for singular knots.  She conjectured that the monoid of singular braids maps injectively into the group algebra of the braid group.  A proof of this conjecture was given by Paris in \cite{Paris}.  Kauffman state models of the Alexander and Jones polynomials were investigated by Fiedler in the context of singular knots \cite{Fiedler}. Quandle theory was extended to the context of non-oriented singular knots in \cite{CEHN1} in order to provide invariants for singular knots.  Further work was carried by the authors in \cite{CCEH} in which we introduced the notion of quandle cocycle invariant for oriented singular knots and links using algebraic structures called \emph{oriented singquandles} and assigning weight functions at both regular and singular crossings.  One of this enhanced invariant's main application was to show that it distinguishes some singular knots with the same number of colorings by an oriented singquanlde, thus making it a \emph{stronger} invariant.

 In \cite{N}, a two-variable polynomial invariant of finite quandles was introduced.  It encodes a set with multiplicities arising from counting trivial actions of elements on other elements of the quandle.  It was further generalized in \cite{N2}.  In \cite{CCE1}, the authors introduced a six-variable polynomial invariant of finite singquandles which extended the quandle polynomial to oriented singquandles.  Furthemore, a \emph{subsingquandle} was introduced and a subsingquandle polynomial was also defined. The polynomial defined for subsingquandles can be thought of as the contributions to the singquandle polynomial coming from the subsingquandles considered.  Using these subsingquandle polynomials, an invariant of singular knots, that generalizes the singquandle counting invariant, was constructed in \cite{CCE1}.

A cohomology theory of quandles was investigated in \cite{CJKLS} and low dimensional cocycles were applied to study state-sum invariants of knots in the $3$-space and knotted surfaces in $4$-space.  In particular, a Boltzmann weight was defined using $3$-cocycles while considering shadow colorings (called also region colorings) in addition to arc colorings thus giving a state sum invariant for knotted surfaces \cite{CJKLS}.  In \cite{CCE} shadow structures for singular knot theory were investigated and used to define invariants of singular knots and links. A notion of an action of a singquandle on a set was introduced to define a shadow counting invariant of singular links which generalize the classical shadow colorings of knots by quandles. Then a shadow polynomial invariant was defined for shadow structures. The shadow counting invariant was enhanced by combining both the shadow counting invariant and the shadow polynomial invariant.

 This article is organized as follows.  In Section~\ref{review}, we review the necessary ingredients of quandles and biquandles and give examples. Section~\ref{OSKQ} surveys 
 oriented singquandles and defines the singquandle counting invariant of singular links.  In Section~\ref{Sec4}, we explain how $2$-cocyles can be obtained from the generalized Reidemeister moves of singular knots and give examples.  This allows us in Section~\ref{inv} to define state-sum invariant of sinuglar knots. Section~\ref{SS1} gives an overview of the construction of the singquandle polynomial, the subsingquandle polynomial, and the polynomial invariant of singular links.  In Section~\ref{SS2}, we survey recent enhancements to the signquandle counting invariant through the use of the singquandle polynomial and the shadow singuqandle polynomial. We end  Section~\ref{Psy} with discussing psyquandle counting invariants of oriented singular links and pseudolinks \cite{NOS} and enhancements  of  the  psyquandle  counting invariant \cite{CN2}.  
 
\section{Review of Quandles}\label{review}
In this section we give a very brief review of quandles and biquandles with some examples.  We start with the following definition of a quandle whose axiomatic come from the three Reidemeister moves on knot diagrams.
\begin{definition}\label{quandledef}
A {\it quandle}, $X$, is a set with a binary operation $(a, b) \mapsto  a * b$ such that

(I) For any $a \in X$,
$a * a =a$.

(II) For any $a,b \in X$, there is a unique $c \in X$ such that
$a= c * b$.

(III)
For any $a,b,c \in X$, we have
$ (a * b) * c=(a * c)* (b * c). $
\end{definition}

\noindent A {\it rack} is a set with a binary operation that
satisfies (II) and (III). Racks and quandles have been studied
extensively in, for example, \cite{Joyce,Matveev}.  For more details on racks and quandles see \cite{EN}.

The following are typical examples of quandles: 
\begin{itemize}

\item
The set $\mathbb{Z}_n$ of integers modulo $n$ is a quandle with the binary operation $x*y=-x+2y$, called dihedral quandle of order $n$. 

\item
Any ${\mathbb{Z} }[t, t^{-1}]$-module $M$ is
a quandle with $a * b=ta+(1-t)b$, for $a,b \in M$, and is called
an {\it  Alexander  quandle}. The previous example corresponds to $M=\mathbb{Z}_n$ and $t=-1$.

\item
A group $G$ with
conjugation as the quandle operation: $a * b = b^{-1} a b$,
denoted by $X=$ Conj$(G)$, is a quandle. 

\item
Any group $G$ with the quandle operation: $a * b = ba^{-1}  b$ is a quandle called Core(G). 
\end{itemize}


Biquandles are generalizations of quandles and are motivated by colorings of \emph{semi-arcs} in a knot or link diagram.  A biquandle is an algebraic structure with two binary operations.
We recall the definition below and provide examples.
\begin{definition}
 A set $(X, \overline{\rhd}, \unrhd)$ is called a {\it biquandle structure} if the following three identities are satisfied.

 (I) For all $x \in X$, $x \unrhd x = x \overline{\rhd} x.$

(II) The maps $\alpha_x, \beta_x: X \to X$ and $S: X \times X \to X \times X$ defined by $\alpha_x(y) = y \overline{\rhd} x,
\beta_x(y) = y \unrhd x $ and $ S(x,y) = (y \overline{\rhd} x, x \unrhd y)$ are invertible, and \label{axiom_2}

(III) The exchange laws satisfied:
\begin{eqnarray*}
 (x \unrhd y) \unrhd (z \unrhd y) &=& (x \unrhd z) \unrhd (y \overline{\rhd} z),\\ 
 (x \unrhd y) \overline{\rhd} (z \unrhd y) &=& (x \overline{\rhd} z) \unrhd (y \overline{\rhd} z),  \\ 
 (x \overline{\rhd} y) \overline{\rhd} (z \overline{\rhd} y) &=& (x \overline{\rhd} z) \overline{\rhd} (y \unrhd z). \label{axiom_3} 
\end{eqnarray*}

\end{definition}

\noindent The following are two typical examples:
\begin{itemize}
\item For any set $X$ and bijection $\sigma: X \to X$ the operations $x \unrhd y = x \overline{\rhd} y = \sigma(x)$ define a biquandle structure called a {\it constant action biquandle.} If $\sigma$ is the identity then we have a {\it trivial biquandle.}

\item Let $X$ be any module over $\mathbb{Z }[t^{\pm 1},s^{\pm 1}].$ Then $X$ is a biquandle with operations $$ x \unrhd y = tx + (s-t)y, \, \,\,  x \overline{\rhd} y =s x $$ known as {\it Alexander biquandle.}
\end{itemize}

\section{Oriented Singular Knots and Singquandles}\label{OSKQ}
\begin{figure}[h] 
\tiny{
\centering
    \includegraphics[scale=0.6]{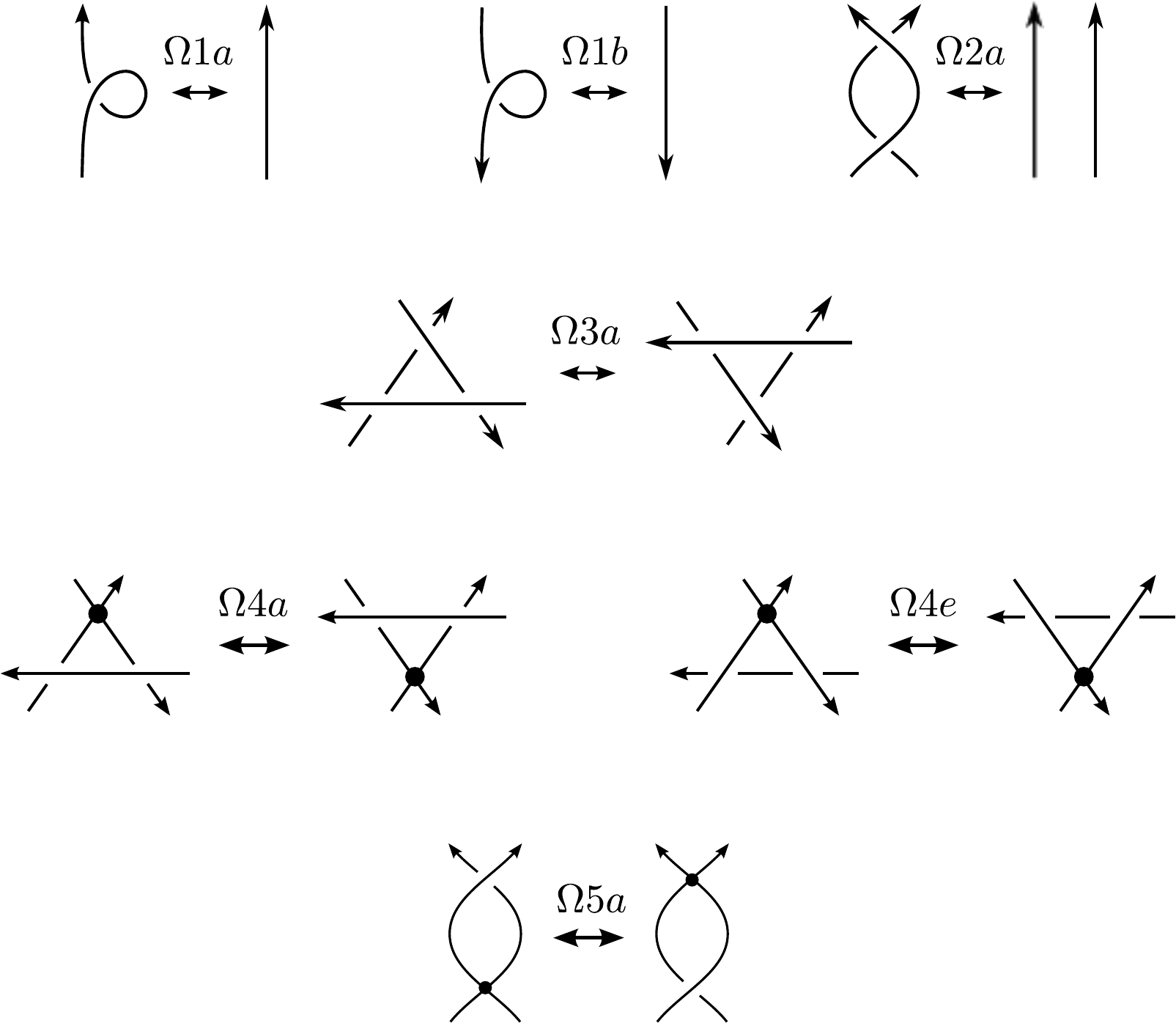}
    \caption{Generating set of singular Reidemeister moves}
    \label{generatingset}}
\end{figure}

Recently, algebraic structures related to singular knot theory were introduced in \cite{CEHN2, BEHY,NOS }. In \cite{BEHY}, generating sets of singular Reidemister moves were introduced with the goal of defining an invariant of singular links. Refer to Figure~\ref{generatingset} for the generating set of the generalized Reidemeister moves used in \cite{BEHY} to define oriented singquandles.


\begin{figure}[h]
	\tiny{
		\centering
		{\includegraphics[scale=0.55]{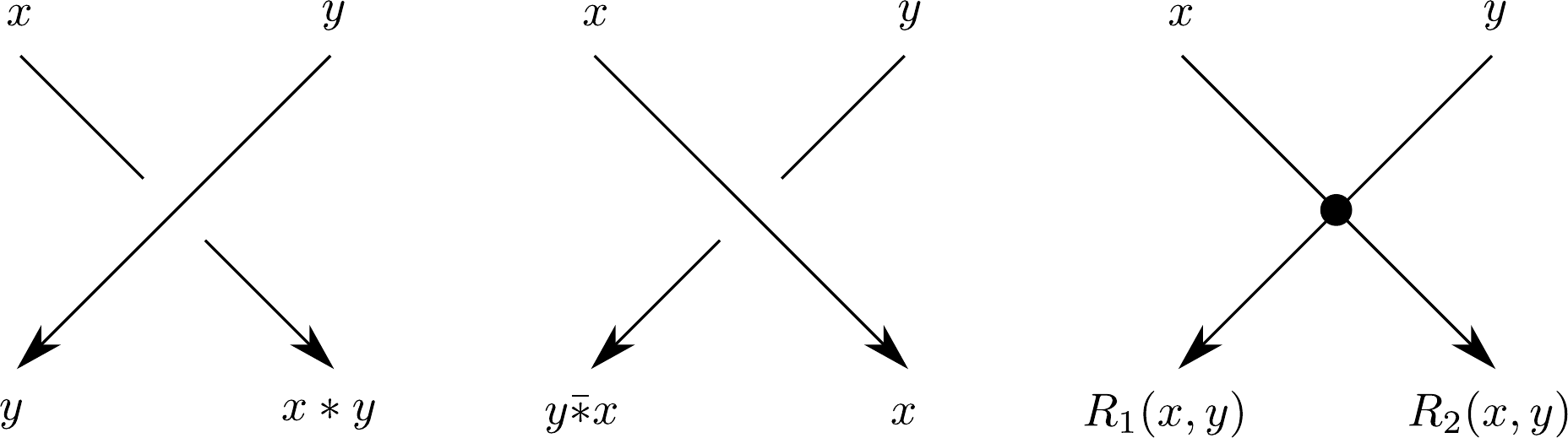}}
		\vspace{.2in}
		\caption{Colorings of classical and singular crossings}
		\label{ColSing}}
\end{figure}

The generating set of Reidemeister moves allows one to introduce four functions as shown in Figure \ref{ColSing} and are used in the definition of an oriented singquandle. The defining axioms of an oriented singquandle were derived from the generalized set of singular Reidemeister moves in Figure~\ref{generatingset} following the crossing rule in Figure \ref{ColSing}.

\begin{definition}\cite{BEHY} \label{oriented SingQdle}
	Let $(X, *)$ be a quandle.  Let $R_1$ and $R_2$ be two maps from $X \times X$ to $X$.  The triple $(X, *, R_1, R_2)$ is called an \emph{oriented singquandle} if the following axioms are satisfied
	
	\begin{eqnarray}
		R_1(x\bar{*}y,z)*y&=&R_1(x,z*y) \ \label{eq1}\\
		R_2(x\bar{*}y, z) & =&  R_2(x,z*y)\bar{*}y \label{eq2}\\
	      (y\bar{*}R_1(x,z))*x   &=& (y*R_2(x,z))\bar{*}z  \label{eq3}\\
R_2(x,y)&=&R_1(y,x*y)  \label{eq4}\\
R_1(x,y)*R_2(x,y)&=&R_2(y,x*y). \label{eq5}	
\end{eqnarray}
\end{definition}

The notion of a homomorphism and isomorphism of oriented singquandles are natural and are given in the following definition.

\begin{definition}\label{singHom}
A map $f: X \rightarrow Y$ is called a \emph{homomorphism} of oriented singquandles $(X, *, R_1, R_2)$ and $(Y, \triangleright, R'_1, R'_2)$ if the following conditions are satisfied for all $x,y \in X$
\begin{eqnarray}
f(x*y)&=&f(x) \triangleright f(y)\label{3.6}\\
f(R_1(x,y))&=&R'_1(f(x),f(y))\label{3.7}\\
f(R_2(x,y))&=&R'_2(f(x),f(y)).\label{3.8}
\end{eqnarray}
An oriented singquandle \emph{isomorphism} is a bijective oriented singquandle homomorphism. We say two oriented singquandles are \emph{isomorphic} if there exists an oriented singquandle isomorphism between them.
\end{definition}
\noindent The following example can be found in \cite{BEHY} and provides a useful family of oriented singquandles over $\mathbb{Z}_n$.

\begin{example}
Let $n$ be a positive integer, let $a$ be an invertible element in $\mathbb{Z}_n$ and let $b,c \in \mathbb{Z}_n$ such that $(1-a)(1 - b-c)=0$.  Then the binary operations $x*y = ax+(1-a)y$, $ R_1(x,y) = bx + cy$ and $R_2(x,y)= acx + [b+ c(1 - a)]y $ make the triple $(\mathbb{Z}_n,*, R_1,R_2)$ into an oriented singquandle.
\end{example}


\noindent Furthermore, the idea of a coloring of a singular link by an oriented singquandle was introduced in \cite{BEHY}. 

\begin{definition}\label{color}
A \emph{coloring} of an oriented singular link $L$ is a function $C : R \rightarrow X$, where $X$ is
a fixed oriented singquandle and $R$ is the set of semiarcs in a fixed diagram of $L$, satisfying the
conditions given in Figure~\ref{ColSing}.
\end{definition}

The set of colorings of $L$ by an oriented singquandle $S$, denoted $Col_S(L)$, was shown to be an invariant of singular links in \cite{BEHY}. From the set of colorings of a singular link by an oriented singquandle we get the \emph{singquandle counting invariant} of the singular link $L$ with respect to the oriented singquandle $S$, denoted $\#Col_S(L)$, by simply computing the cardinality of $Col_S(L)$ (see \cite{BEHY} for more details).
We will now explore recent enhancements to the signquandle counting invariant, which futher extract information from the set of colorings of a singular link by an oriented singquandle. 
\section{Two Dimensional Cocycles from Singular Knots}\label{Sec4}

Quandle $2$-cocycles proved to be very powerful in constructing state-sum invariant of classical knots and knotted surfaces \cite{CJKLS}.  Similar ideas were developed by the authors in \cite{CCEH}.  The following will be an overview of the results obtained through the use of the state-sum invariant defined using oriented singquandles. For a detailed construction of the state-sum invariant, see \cite{CCEH}. 

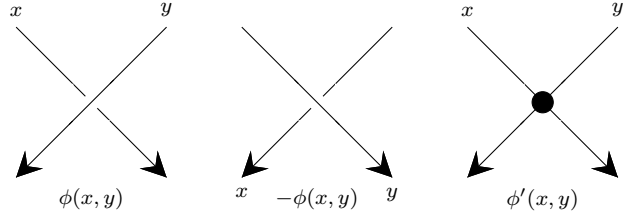
\begin{figure}[h]
\begin{tikzpicture}[use Hobby shortcut]
\begin{knot}[
consider self intersections=true,
  ignore endpoint intersections=false,
  flip crossing/.list={3,4,5,6,11,13},
  clip width = 4
]
\strand[decoration={markings,mark=at position 1 with
    {\arrow[scale=3,>=stealth]{>}}},postaction={decorate}] (1,1) ..(-1,-1); 
\strand[decoration={markings,mark=at position 1 with
    {\arrow[scale=3,>=stealth]{>}}},postaction={decorate}] (-1,1) ..(1,-1);
\strand[decoration={markings,mark=at position 1 with
    {\arrow[scale=3,>=stealth]{>}}},postaction={decorate}] (2,1) ..(4,-1);
\strand[decoration={markings,mark=at position 1 with
    {\arrow[scale=3,>=stealth]{>}}},postaction={decorate}] (4,1) ..(2,-1);  
\strand[decoration={markings,mark=at position 1 with
    {\arrow[scale=3,>=stealth]{>}}},postaction={decorate}] (5,1) ..(7,-1);
\strand[decoration={markings,mark=at position 1 with
    {\arrow[scale=3,>=stealth]{>}}},postaction={decorate}] (7,1) ..(5,-1);      
\end{knot}
\node[above] at (-1,1) {\tiny $x$};
\node[above] at (1,1) {\tiny $y$};
\node[below] at (0,-1) {\tiny $\phi(x,y)$};

\node[below] at (2,-1) {\tiny $x$};
\node[below] at (4,-1) {\tiny $y$};
\node[below] at (3,-1) {\tiny $-\phi(x,y)$};

\node[circle,draw=black, fill=black, inner sep=0pt,minimum size=8pt] (a) at (6,0) {};
\node[above] at (5,1) {\tiny $x$};
\node[above] at (7,1) {\tiny $y$};
\node[below] at (6,-1) {\tiny $\phi'(x,y)$};
\end{tikzpicture}
\vspace{.2in}
		\caption{Boltzmann weights at classical and singular crossings.}
		\label{weights}
\end{figure}

Let $(X,*)$  be a quandle and  $A$ be an abelian group. At crossings (first and second) of Figure \ref{weights} we assign a Boltzmann weight $\phi(x,y)$ to the positive crossing and the weight $-\phi(x,y)$ to the negative crossing.  One then obtains a \textit{$2$-cocycle} of $X$ on $A$ as a function $\phi: X \times X \rightarrow A$ that satisfies some conditions coming from Reidemeister moves.
Precisely, the function $\phi$ satisfies the so called \emph{$2$-cocycle condition} for all $x,y,z \in X,$
\begin{eqnarray*}
\phi(x,y) + \phi (x*y,z)=\phi(x,z) + \phi(x*z,y*z). 
\end{eqnarray*}
This condition is derived from Reidemeister move III. Moreover, Reidemeister move I implies the condition $\phi(x,x)=0$ on a $2$-cocycle. Finally, Reidemeister move II imposes that the total Boltzmann weight should be zero; this condition is automatically satisfied due to the weight at a positive and negative crossing canceling each other (see Figure ~\ref{weights}). 

Naturally, one wants to extend this definition to singular knots and links. To this end, we define a function $\phi^{\prime} :X \times X \to A $ that represents the Boltzmann weight at a singular crossing (see the third diagram in Figure \ref{weights}). 

\begin{figure}[h]
	\tiny{
		\centering
		{\includegraphics[scale=0.55]{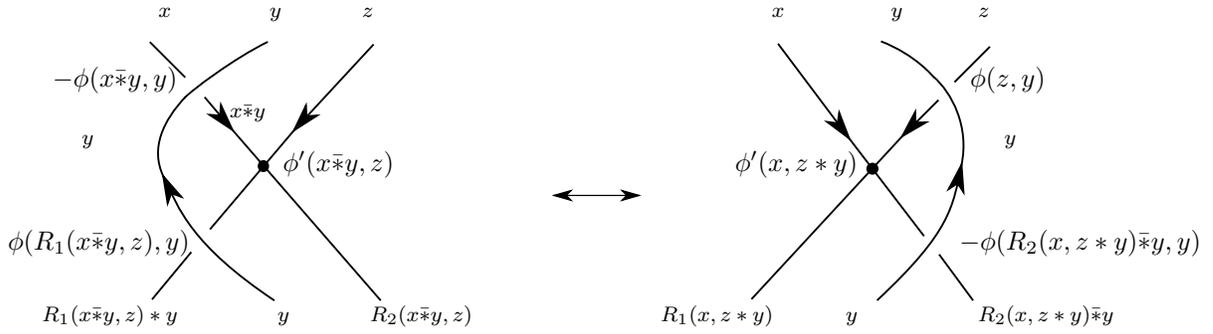}
			\put(-315,108){$x$}
			\put(-273,108){$y$}
			\put(-288,70){$x\bar{*}y$}
			\put(-344,60){$y$}
			\put(5,60){$y$}
			\put(-359,-7){$R_1(x\bar{*}y,z)*y$}
			\put(-270,-7){$y$}
			\put(-235,-7){$R_2(x\bar{*}y,z)$}
			\put(-125,-7){$R_1(x,z*y)$}
			\put(-55,-7){$y$}
			\put(-5,-7){$R_2(x,z*y)\bar{*}y$}
			\put(-238,108){$z$}
			\put(-83,108){$x$}
			\put(-38,108){$y$}
			\put(-5,108){$z$}
            \put(-355,82){\footnotesize{$-\phi(x\bar{*}y,y)$}}
            \put(-8,82){\footnotesize{$\phi(z,y)$}}
           \put(-97,50){\footnotesize{$\phi^{\prime}(x,z*y)$}}
             \put(-268,50){\footnotesize{$\phi^{\prime}(x\bar{*}y,z)$}}
            \put(-12,20){\footnotesize{$-\phi(R_2(x,z*y)\bar{*}y,y)$}}
            \put(-372,20){\footnotesize{$\phi(R_1(x\bar{*}y,z),y)$}}
		}
		\vspace{.2in}
		\caption{The Boltzmann weights and Reidemeister move $\Omega 4a$.} 
		\label{W4a}}
\end{figure}
The two functions $\phi$ and $\phi'$ satisfy the conditions under singular Reidemeister moves. Using Figure \ref{W4a} we obtain the following condition:
\begin{equation}\label{troisieme}
-\phi(x\bar{*}y,y)+\phi^{\prime}(x\bar{*}y,z)+\phi(R_1(x\bar{*}y,z),y)=\phi(z,y)+\phi^{\prime}(x,z*y)-\phi(R_2(x,z*y)\bar{*}y,y).
\end{equation}

\begin{figure}[h]
	\tiny{
		\centering
		{\includegraphics[scale=0.55]{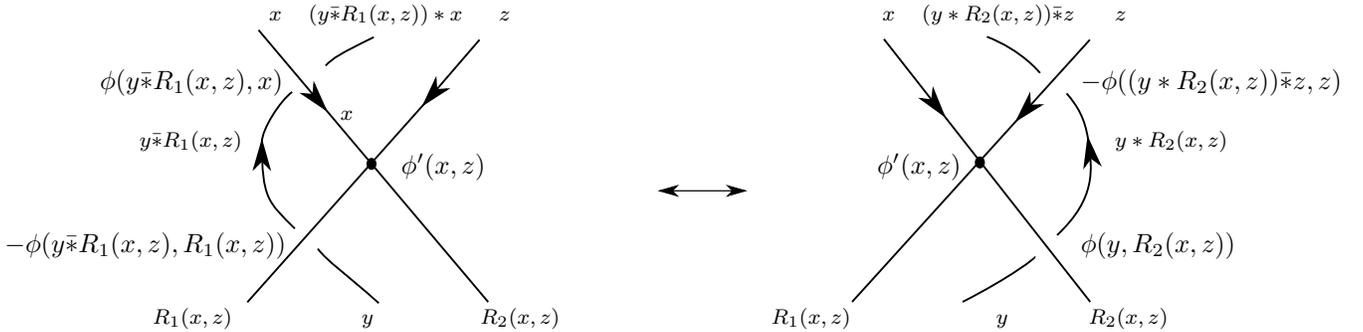}
			\put(-315,108){$x$}
			\put(-300,108){$(y\bar{*}R_1(x,z))*x$}
			\put(-288,70){$x$}
            \put(-379,82){\footnotesize{$\phi( y\bar{*}R_1(x,z),x)$}}
            \put(-8,82){\footnotesize{$-\phi((y*R_2(x,z))\bar{*}z, z)$}}
            \put(-85,50){\footnotesize{$\phi^{\prime}(x,z)$}}
             \put(-265,50){\footnotesize{$\phi^{\prime}(x,z)$}}
            \put(-8,20){\footnotesize{$\phi(y,R_2(x,z))$}}
            \put(-415,20){\footnotesize{$-\phi(y\bar{*}R_1(x,z),R_1(x,z))$}}
			\put(-364,60){$y\bar{*}R_1(x,z)$}
			\put(5,60){$y*R_2(x,z)$}
			\put(-359,-7){$R_1(x,z)$}
			\put(-280,-7){$y$}
			\put(-235,-7){$R_2(x,z)$}
			\put(-125,-7){$R_1(x,z)$}
			\put(-40,-7){$y$}
			\put(-5,-7){$R_2(x,z)$}
			\put(-228,108){$z$}
			\put(-83,108){$x$}
			\put(-68,108){$(y*R_2(x,z))\bar{*}z$}
			\put(5,108){$z$}
		}
		\vspace{.2in}
		\caption{The Boltzmann weights and Reidemeister move $\Omega 4e$. }
		\label{W4b}}
\end{figure}
Figure \ref{W4b} implies the following equation:
\begin{eqnarray}\label{RR12}
\phi(y\bar{*}R_1(x,z),x)-\phi(y\bar{*}R_1(x,z),R_1(x,z))=
-\phi((y*R_2(x,z))\bar{*}z, z)+\phi(y,R_2(x,z)),
\end{eqnarray}

\begin{figure}[h]
\tiny{
  \centering
   {\includegraphics[scale=0.39]
   {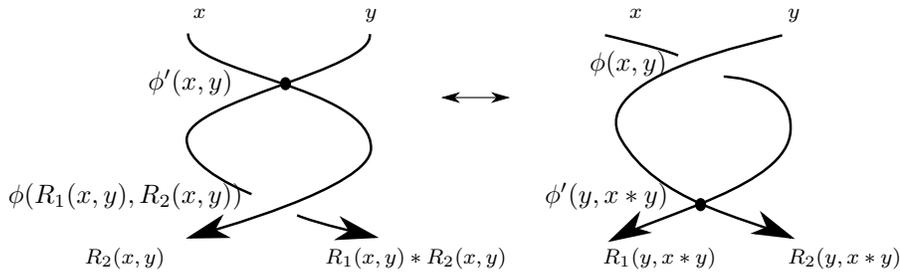}
	\put(-163,85){$y$}
    \put(-228,85){$x$}
    \put(-245,58){\footnotesize{$\phi^{\prime}(x,y)$}}
    \put(-298,15){\footnotesize{$\phi(R_1(x,y),R_2(x,y))$}}
    \put(-78,65){\footnotesize{$\phi(x,y)$}}
    \put(-95,15){\footnotesize{$\phi^{\prime}(y,x*y)$}}
	\put(-3,-8){$R_2(y,x*y)$}
	 \put(-269,-8){$R_2(x,y)$}  
    \put(-63,85){$x$} 
    \put(-3,85){$y$}   
    \put(-73,-8){$R_1(y,x*y)$}
    \put(-178,-8){$R_1(x,y)*R_2(x,y)$}}
     \vspace{.2in}
     \caption{The Reidemeister move $\Omega 5a$ and colorings}
     \label{W5}}
\end{figure}
while Figure \ref{W5} implies the following equation:
\begin{equation}\label{deuxiemeequation}
\phi^{\prime}(x,y)+\phi(R_1(x,y),R_2(x,y))=\phi(x,y)+\phi^{\prime}(y,x*y).
\end{equation}

\section{Cocycle Invariants of Singular Knots and Links}\label{inv}


The definition of the cocycle invariant was extended by the authors to singular knots \cite{CCEH} where it was proved that it is indeed an invariant for singular knots.  We start by recalling that given an abelian group $A$ and a quandle $(X,*)$, a $2$-cocycle is a function $\phi: X \times X \rightarrow A$ that satisfies the following for all $x,y,z \in X$
\begin{eqnarray}\label{standard}
\phi(x,y) + \phi (x*y,z)&=&\phi(x,z) + \phi(x*z,y*z) \text{ and } \phi(x,x)=0.
\end{eqnarray}
In order to construct our invariant we need to solve the system made of  equations~(\ref{troisieme}),\\ ~(\ref{RR12}), ~(\ref{deuxiemeequation}) and~(\ref{standard}).  First we need to recall the notion of coloring of a knot diagram by an oriented singquandle. Let $D$ be a diagram of a singular knot and let $(X, *, R_1, R_2)$ be an oriented singquandle, then a coloring of $D$ by $(X, *, R_1, R_2)$ is defined in a similar way to the case of colorings of classical knots by quandles.  To be precise, the colorings at positive, negative and singular crossings are given by Figure~\ref{ColSing}.
As in the case of classical knots we will assign the weights $\phi(x,y)$ and $-\phi(x,y)$ at a positive and negative crossing respectively as shown Figure~\ref{weights}. Furthermore, we associate the weight $\phi'(x,y)$ at a singular crossing as shown in Figure~\ref{weights}.  Now, given an abelian group $A$ (denoted multiplicatively) and cocycles $\phi, \phi': X \times X \rightarrow A$, we define the state sum in exactly the same manner as state sum for classical knots (see \cite{CJKLS}
page 3953). It will be clear below that the state sum does not depend on the choice of a diagram of a given knot or link.
\begin{definition}\label{StSum}\cite{CCEH}
Let $(X, *, R_1, R_2)$ be an oriented singquandle and let $A$ be an abelian group.  Assume that $\phi, \phi': X \times X \rightarrow A$ satisfy the  equations~(\ref{troisieme}),~(\ref{RR12}), ~(\ref{deuxiemeequation}) and~(\ref{standard}).  Then the state sum for a diagram of a singular knot $K$ is given by $\displaystyle \Phi_{\phi, \phi'}(K)=\sum_{\mathcal{C}}\prod_{\tau} \psi(x,y) $, where the product is taken over all classical and singular crossings $\tau$ of the digaram of $K$, the sum is taken over all possible colorings $\mathcal{C}$ of the knot $K$.   
\end{definition}
Notice that in this definition $\psi(x,y)= 
    \phi(x,y)^{\pm 1}$ at classical positive or negative crossing and $\psi(x,y)= \phi'(x,y)$ at a singular crossing as in Figure~\ref{weights}.
  The following theorem stating that $\Phi_{\phi, \phi'}(K)$ is an invariant of singular knots was proved in \cite{CCEH}.   
\begin{theorem}
Let $\phi, \phi': X \times X \rightarrow A$ be maps satisfying the conditions of Definition~\ref{StSum}.  The \emph{state sum} associated with $\phi$ and $\phi'$ is invariant under the moves listed in the \emph{generating set} of singular Reidemester moves in Figure~\ref{generatingset}, so it defines an invariant of singular knots and links.
\end{theorem}



\noindent We provide an example of singquandles and their $2$-cocycles. We then use the cocycle invariant introduced in Definition~\ref{StSum} written in the form $\displaystyle \Phi_{\phi, \phi'}(K)=\sum_{i=0}^{n-1}a_iu^i$ to distinguish a pair of singular knots and links.  For more examples the reader can consult \cite{CCEH}. 

\begin{example}
Let $(\mathbb{Z}_{6}, *, R_1, R_2)$ be the oriented singquandle with $x *y = -x+2y = x \bar{*}y$ and $R_1(x,y) =3+2x-y$ and $R_2(x,y) = 3$ and weight functions $\phi(x,y) = 2x+3x^2 -2y-xy -2y^2$ and $\phi'(x,y) = 3+x+x^2+2y-xy$. In \cite{CCEH}, the state-sum invariant was used to distinguish the following singular knots listed as $5^k_6$ and $5^k_7$ in \cite{Oyamaguchi}. 

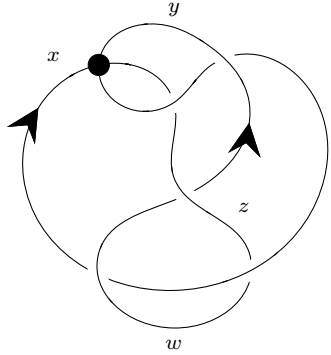
\begin{figure}[h]
\begin{tikzpicture}[use Hobby shortcut]
\begin{knot}[thick,
consider self intersections=true,
  ignore endpoint intersections=false,
  flip crossing/.list={3,4,5,6,11,13}
]
\strand ([closed]0,2) [decoration={markings,mark=at position .5 with
    {\arrow[scale=3,>=stealth]{>}}},postaction={decorate}]..(-1,1.5)..(0,1)..(.5,1.5)..(2,0)..(0,-1.5)..(-2,0)[decoration={markings,mark=at position .93 with
    {\arrow[scale=3,>=stealth]{>}}},postaction={decorate}]..(-1,1.5)..(0,1)..(0,0)..(1,-1)..(0,-2)..(-1,-1)..(0,-.3)..(1,1)..(.4,1.8)..(0,2);
\end{knot}

\node[circle,draw=black, fill=black, inner sep=0pt,minimum size=8pt] (a) at (-1,1.5) {};
\node[above] at (-1.6,1.4) {\tiny $x$};
\node[above] at (0,2) {\tiny $y$};
\node[below] at (0,-2) {\tiny $w$};
\node[right] at (.7,-.4) {\tiny $z$};
\end{tikzpicture}
\vspace{.2in}
		\caption{Diagram for $5^k_6$.}
		\label{5k6}
\end{figure}

\begin{figure}[h]

\begin{tikzpicture}[use Hobby shortcut]
\begin{knot}[thick,
consider self intersections=true,
  ignore endpoint intersections=false,
  flip crossing/.list={3,4,5,6,11,13}
]
\strand ([closed]0,2)[decoration={markings,mark=at position .45 with
    {\arrow[scale=3,>=stealth]{<}}},postaction={decorate}]..(-1,1.5)..(0,1)..(.5,1.5)..(2,0)..(0,-1.5)..(-2,0)[decoration={markings,mark=at position .93 with
    {\arrow[scale=3,>=stealth]{<}}},postaction={decorate}]..(-1,1.5)..(0,1)..(0,0)..(1,-1)..(0,-2)..(-1,-1)..(0,-.3)..(1,1)..(.4,1.8)..(0,2);
\end{knot}

\node[circle,draw=black, fill=black, inner sep=0pt,minimum size=8pt] (a) at (.65,1.55) {};
\node[above] at (0,2) {\tiny $x$};
\node[above] at (1.2,1.6) {\tiny $y$};
\node[right] at (.7,-.4) {\tiny $w$};
\node[left] at (-2,0) {\tiny $z$};
\end{tikzpicture}
\vspace{.2in}
		\caption{Diagram for $5^k_7$.}
		\label{5k7}
\end{figure}
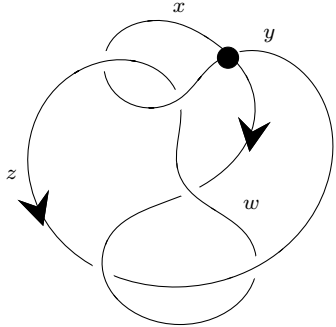

The two singular knots have 6 colorings by this given oriented singquandle.  On the other hand, the state-sum invariant can be used to distinguish them, $\Phi_{\phi,\phi'}(5^k_6) = 6u^3 \neq 6= \Phi_{\phi,\phi'}(5^k_7)$.
\end{example}

\section{Singquandle Polynomial and the Singular Link Invariant }\label{SS1}

 The authors of \cite{CCE1} took a different approach to extract additional information from the singquandle coloring set. In \cite{CCE1}, the authors defined the singquandle polynomial, the subsingquandle polynomial and a polynomial invariant of singular links. These polynomials were used to study how subsingquadles of a singquandle sit in the singquandle. In this section, we will give an overview of the construction of the singquandle polynomial, the subsingquandle polynomial, and the polynomial invariant of singular links. We will follow the construction and notation introduced in \cite{CCE1}.

\begin{definition}
Let $(X,*,R_1,R_2)$ be a finite singquandle. For every $x \in X$, define the following sets
\[ C^1(x) = \lbrace y \in X \, \vert \, y * x = y \rbrace \quad \text{and} \quad R^1(x) = \lbrace y \in X \, \vert \, x * y = x \rbrace, \]
\[ C^2(x) = \lbrace y \in X \, \vert \, R_1(y , x) = y \rbrace \quad \text{and} \quad R^2(x) = \lbrace y \in X \, \vert \, R_1(x , y) = x \rbrace, \]
\[ C^3(x) = \lbrace y \in X \, \vert \, R_2(y , x) = y \rbrace \quad \text{and} \quad R^3(x) = \lbrace y \in X \, \vert \, R_2(x , y) = x \rbrace. \]\\
Let $c^i(x) = \vert C^i(x)\vert$ and $r^i(x) = \vert R^i(x)\vert$ for $i=1,2,3$. Then the \emph{singquandle polynomial of X} is defined by

\[ sqp(X) = \sum_{x\in X} s_1^{r^1(x)}t_1^{c^1(x)}s_2^{r^2(x)}t_2^{c^2(x)}s_3^{r^3(x)}t_3^{c^3(x)}.  \]
\end{definition}
We note that the value $r^i(x)$ is the number of elements in $X$ that act trivially on $x$, while $c^i(x)$ is the number of elements of $X$ on which $x$ acts trivially via $*, R_1$ and $R_2.$ Furthermore, if $Y \subset X$ is a subsingquandle we can define the following singquandle polynomial for $Y$ as a subsignquandle of $X$.

{\color{red} }
\begin{definition}
Let $(X,*, R_1,R_2)$ be a finite singquandle and $S \subset  X$ a subsingquandle. Then the \emph{subsingquandle polynomial} is 
\[ Ssqp(S \subset X ) = \sum_{x \in S} s_1^{r^1(x)}t_1^{c^1(x)}s_2^{r^2(x)}t_2^{c^2(x)}s_3^{r^3(x)}t_3^{c^3(x)}. \]
\end{definition} 
The subsingquandle polynomial can be thought of as the contribution to the singquandle polynomial coming from the subsingquandle we are considering. Using the subsingquandle polynomial it was shown that the following polynomial is an invaraint of singular links in \cite{CCE1}.

\begin{definition}
Let $L$ be a singular link, $(X,*,R_1,R_2)$ a finite singquandle. Then the multiset
\[ \Phi_{Ssqp}(L,X) = \lbrace Ssqp(Im(f) \subset X) \, \vert \, f \in \text{Hom}(\mathcal{SQ}(L), X\rbrace \]
is the \emph{subsingquandle polynomial invariant of $L$} with respect to $X$. We can also represent this invariant in the following polynomial-style form by converting the multiset elements to exponents of a formal variable $u$ and converting their multiplicities to coefficients:
\[ \phi_{Ssqp}(L,X) = \sum_{f \in \textup{Hom}(\mathcal{SQ}(L),X)} u^{Ssqp(Im(f) \subset X)}.\]
\end{definition}

The authors of \cite{CCE1} used the subsingquandle polynomial invariant to distinguish several pairs of singular knots with the same singquandle counting invariant. The following example and several others examples with explicit computations can be found in \cite{CCE1}.
\begin{example}
Consider the singquandle $X = \mathbb{Z}_8$ with operations $x*y = 5x+4y = x\bar{*} y$ and $R_1(x,y)= 6+5x+6xy$, and $R_2(x,y)= 6+5y+6xy$. The knot $K_1$, Figure~\ref{k1}, and the knot $K_2$, Figure~\ref{k2}, both have 8 colorings by the given oriented singquandle.

\begin{figure}[ht]
\centering
\begin{tikzpicture}[use Hobby shortcut,scale=.7]
\begin{knot}[
  consider self intersections=true,
  ignore endpoint intersections=false,
  flip crossing=3,
  only when rendering/.style={
  }
  ]
\strand ([closed]0,1.5)[decoration={markings,mark=at position .5 with
    {\arrow[scale=3,>=stealth]{<}}},postaction={decorate}]..(-1.2,-1.5).. (3,-3.5) ..(0,-1.2) ..(-3,-3.5) ..(1.2,-1.5)..(0,1.5);
\end{knot}
\node[circle,draw=black, fill=black, inner sep=0pt,minimum size=8pt] (a) at (0,-3.65) {};
\node[circle,draw=black, fill=black, inner sep=0pt,minimum size=8pt] (a) at (1.25,-1.3) {};
\node[circle,draw=black, fill=black, inner sep=0pt,minimum size=8pt] (a) at (-1.25,-1.3) {};
\node[above] at (0,-1) {\tiny $x$};
\node[above] at (0,1.5) {\tiny $y$};
\node[right] at (1,-2) {\tiny $z$};
\node[right] at (3,-3.5) {\tiny $w$};
\node[left] at (-3,-3.5) {\tiny $l$};
\node[left] at (-1,-2) {\tiny $k$};

\end{tikzpicture}
		\caption{Diagram of $K_1$.}
		\label{k1}
\end{figure}
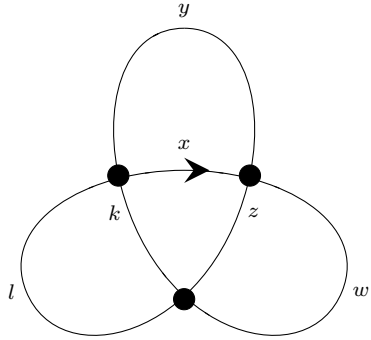

\begin{figure}[ht]
\centering
\begin{tikzpicture}[use Hobby shortcut,scale=.7]
\begin{knot}[
  consider self intersections=true,
  ignore endpoint intersections=false,
  flip crossing=4,
  clip width=5,
  only when rendering/.style={
  }
  ]
\strand ([closed]0,1.5)[decoration={markings,mark=at position .5 with
    {\arrow[scale=3,>=stealth]{<}}},postaction={decorate}]..(-1.2,-1.5).. (3,-3.5) ..(0,-1.2) ..(-3,-3.5) ..(1.2,-1.5)..(0,1.5);
\end{knot}
\node[circle,draw=black, fill=black, inner sep=0pt,minimum size=9.5pt] (a) at (1.25,-1.35) {};
\node[circle,draw=black, fill=black, inner sep=0pt,minimum size=9.6pt] (a) at (-1.241,-1.35) {};
\node[above] at (0,-1) {\tiny $x$};
\node[above] at (0,1.5) {\tiny $y$};
\node[below] at (-1.5,-4.5) {\tiny $z$};
\node[right] at (3,-3.5) {\tiny $w$};
\node[left] at (-1,-2.5) {\tiny $k$};
\end{tikzpicture}
		\caption{Diagram of $K_2$.}
		\label{k2}
\end{figure}
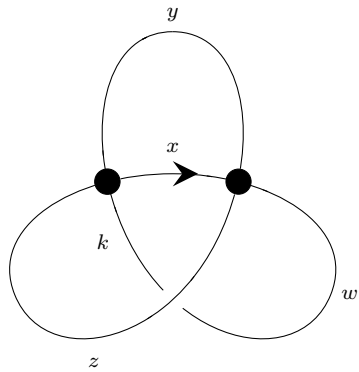
However, the subsingquandle polynomial  invariants
\[ \phi_{Ssqp}(K_1,X) = 4 u^{s_1^4 t_1^4 s_2^2 t_2^2 s_3 t_3}+4 u^{2 s_1^4 t_1^4 s_2^2 t_2^2 s_3
   t_3}\]
and
\[ \phi_{Ssqp}(K_2,X)=4 u^{4 s_1^4 t_1^4 s_3 t_3}+4 u^{s_1^4 t_1^4 s_2^2 t_2^2 s_3 t_3}\]
distinguish the singular knots $K_1$ and $K_2$.
\end{example}

\section{Singquandle Shadows, Singquandle Shadow Polynomial and the Singular Link Invariant }\label{SS2}
\begin{figure}[h]
\begin{tikzpicture}[use Hobby shortcut,scale=1.5]
\begin{knot}[
consider self intersections=true,
clip width=5,
  ignore endpoint intersections=false,
  flip crossing/.list={3,4,5,6,11,13}
]
\strand[decoration={markings,mark=at position 1 with
    {\arrow[scale=1,>=stealth]{>}}},postaction={decorate}] (1,1) ..(-1,-1.1); 
\strand[decoration={markings,mark=at position 1 with
    {\arrow[scale=1,>=stealth]{>}}},postaction={decorate}] (-1,1) ..(1,-1.1);
    
\strand[decoration={markings,mark=at position 1 with
    {\arrow[scale=1,>=stealth]{>}}},postaction={decorate}] (3,1) ..(5,-1.1);
\strand[decoration={markings,mark=at position 1 with
    {\arrow[scale=1,>=stealth]{>}}},postaction={decorate}] (5,1) ..(3,-1.1);

\draw[decoration={markings,mark=at position 1 with
    {\arrow[scale=1,>=stealth]{>}}},postaction={decorate}] (1,-2) ..(3,-4.1);
\draw[decoration={markings,mark=at position 1 with
    {\arrow[scale=1,>=stealth]{>}}},postaction={decorate}] (3,-2) ..(1,-4.1);     
\end{knot}

\node[above] at (-1,1) {\tiny $s_1$};
\node[above] at (1,1) {\tiny $s_2$};

\node[below] at (-1,-1.1) {\tiny $s_2$};
\node[below] at (1,-1.1) {\tiny $s_1 *s_2$};

\draw (-.5,0) node [draw] {\tiny $x$};
\draw (0,.6) node [draw] {\tiny $x \cdot s_1$};
\draw (1.25,.3) node [draw] {\tiny $(x \cdot s_1) \cdot s_2$};
\draw  (0,-.6) node [draw] {\tiny $x \cdot s_2$};
\draw (1.5,-.3) node [draw] {\tiny $(x \cdot s_2) \cdot (s_1 * s_2)$};
\node at (1,0) {\tiny $=$};
\node[below] at (3,-1.1) {\tiny $s_1$};
\node[below] at (5,-1.1) {\tiny $s_2$};

\node[above] at (3,1) {\tiny $s_2$};
\node[above] at (5,1) {\tiny $s_1 * s_2$};

\draw (3.5,0) node [draw] {\tiny $x$};
\draw (4,.6) node [draw] {\tiny $x \cdot s_2$};
\draw (5.5,.3) node [draw] {\tiny $(x \cdot s_2) \cdot (s_1 * s_2)$};
\draw (4,-.6) node [draw] {\tiny $x \cdot s_1$};
\draw (5.2,-.3) node [draw] {\tiny $(x \cdot s_1) \cdot s_2$};
\node at (5,0) {\tiny $=$};

\node[above] at (1,-2) {\tiny $s_1$};
\node[above] at (3,-2) {\tiny $s_2$};

\node[below] at (1,-4.1) {\tiny $R_1(s_1,s_2)$};
\node[below] at (3,-4.1) {\tiny $R_2(s_1,s_2)$};

\draw (1.5,-3) node [draw] {\tiny $x$};
\draw (2,-2.4) node [draw] {\tiny $x \cdot s_1$};
\draw (3.3,-2.7) node [draw] {\tiny $(x \cdot s_1) \cdot s_2$};
\draw (2,-3.97) node [draw] {\tiny $x \cdot R_1(s_1,s_2)$};
\draw (4,-3.3) node [draw] {\tiny $(x \cdot R_1(s_1,s_2)) \cdot R_2(s_1,s_2) $};
\node at (3,-3) {\tiny $=$};
\node[circle,draw=black, fill=black, inner sep=0pt,minimum size=6pt] (a) at (2,-3.05) {};
\end{tikzpicture}
	\caption{Shadow coloring at positive, negative, and singular crossings.}
		\label{shadowX}
\end{figure}
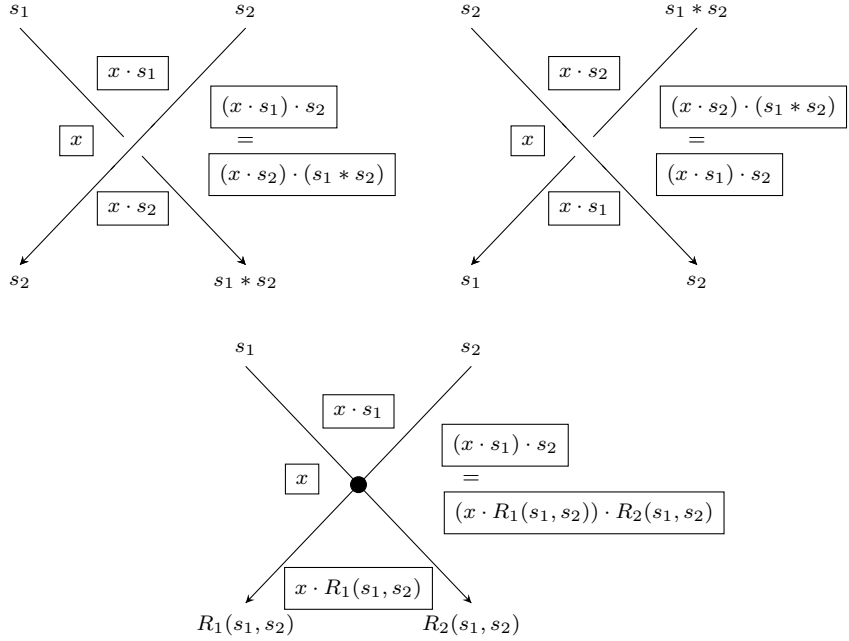

 The ideas introduced in \cite{CCE1} were further studied and generalized in \cite{CCE}.
Precisely, in \cite{CCE}, shadow structures for singular knot theory were investigated and used to define invariants of singular knots and links. A notion of an action of a singquandle on a set was introduced to define a shadow counting invariant of singular links which generalize the classical shadow colorings of knots by quandles. Then a shadow polynomial invariant was defined for shadow structures. The shadow counting invariant was enhanced by combining both the shadow counting invariant and the shadow polynomial invariant.  The following is a brief summary of the results from \cite{CCE}.
\begin{definition}\cite{CCE}
Let $(S,*,R_1,R_2)$ be a singquandle. An \emph{S-set} is a set  $X$ and a map $\cdot : X \times S \rightarrow X$ satisfying the following conditions: 

   (I) For all $s \in S$, $\cdot s : X \rightarrow X$ mapping $x$ to $x \cdot s$ is a bijection.
   
    (II) For all $s_1, s_2 \in S$ and $x \in X$,
 \begin{eqnarray}
(x \cdot s_1)\cdot s_2 &=& (x \cdot s_2)\cdot ( s_1 * s_2)\\
(x \cdot s_1) \cdot s_2 &=& (x \cdot R_1(s_1,s_2)) \cdot R_2(s_1,s_2).
\end{eqnarray}

\end{definition}
The meaning of these two equations will become clear from Figure~\ref{shadowX}.

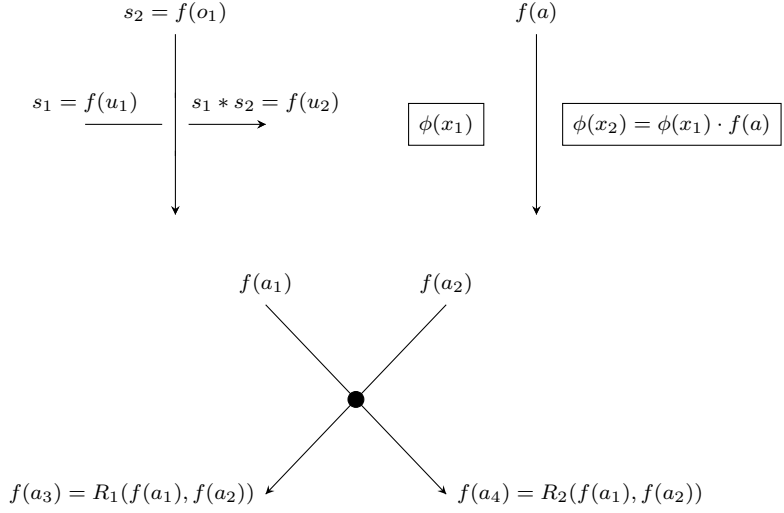
\begin{figure}[h]
\begin{tikzpicture}[use Hobby shortcut,scale=1.2]
\begin{knot}[
consider self intersections=true,
clip width=5,
  ignore endpoint intersections=false,
  flip crossing/.list={3,4,5,6,11,13}
]
\strand[decoration={markings,mark=at position 1 with
    {\arrow[scale=1,>=stealth]{>}}},postaction={decorate}] (0,1) ..(0,-1); 
\strand[decoration={markings,mark=at position 1 with
    {\arrow[scale=1,>=stealth]{>}}},postaction={decorate}] (-1,0) ..(1,0);
    
\strand[decoration={markings,mark=at position 1 with
    {\arrow[scale=1,>=stealth]{>}}},postaction={decorate}] (4,1) ..(4,-1);

\draw (1,-2)..(3,-4.1)[decoration={markings,mark=at position 1 with
    {\arrow[scale=1,>=stealth]{>}}},postaction={decorate}];
\draw[decoration={markings,mark=at position 1 with
    {\arrow[scale=1,>=stealth]{>}}},postaction={decorate}] (3,-2)..(1,-4.1);
   
\end{knot}

\node[above] at (-1,0) {\tiny $s_1=f(u_1)$};
\node[above] at (0,1) {\tiny $s_2=f(o_1)$};
\node[above] at (1,0) {\tiny $s_1 * s_2 =f(u_2)$};

\draw (3,0) node [draw] {\tiny $\phi(x_1)$};
\draw (5.5,0) node [draw] {\tiny $\phi(x_2) = \phi(x_1)\cdot f(a)$};
\node[above] at (4,1) {\tiny $f(a)$};

\node[above] at (1,-2) {\tiny $f(a_1)$};
\node[above] at (3,-2) {\tiny $f(a_2)$};
\node[right] at (3,-4.1) {\tiny $f(a_4)=R_2 (f(a_1),f(a_2))$};
\node[left] at (1,-4.1) {\tiny $f(a_3)=R_1(f(a_1),f(a_2))$};
\node[circle,draw=black, fill=black, inner sep=0pt,minimum size=6pt] (a) at (2,-3.05) {};
\end{tikzpicture}
	\caption{Arcs and regions of diagram $D$.}
		\label{shadowrule}
\end{figure}

\begin{definition}\cite{CCE}
A \emph{singquandle shadow} or $S$-shadow is the pair of an oriented sinquandle $(S,*,R_1,R_2)$ and a $S$-set $(X,\cdot)$, denoted by $(S,X,*,R_1,R_2,\cdot)$ or simply by $(S,X)$. Let $S^{\prime}$ be a subsingquandle of $S$. A subset $Y$ of $X$ closed under the action of $S^{\prime}$ is an \emph{subshadow} of $(S,X)$, which we will denote by $(S^\prime, Y) \subset (S,X)$.
\end{definition}


 

As in the case of classical knot theory, the authors of \cite{CCE1} introduced the idea of a fundamental singquandle associated with a singular link to define the notion of colorings of that singular link by a singquandle.  Thus leading to invariants of singular links.
We give the definition of a homomorphism of singquandle shadows.

\begin{definition}\cite{CCE}
A \emph{homomorphism of sinquandle shadows} between $(S,X,*,R_1,R_2,\cdot)$ and \\ $(S',Y,\triangleright,R'_1,R'_2,\bullet)$ is a pair of maps $\phi:(X, \cdot) \rightarrow (Y,\bullet)$ and $f:(S,*,R_1,R_2) \rightarrow (S',\triangleright,R'_1,R'_2)$, such that $f$ is a singquandle homomorphism, that is the identities (\ref{3.6}), (\ref{3.7}) and (\ref{3.8}) are satisfied and for all $x \in X$ and $s \in S$, we have
\begin{eqnarray}\label{action}
\phi(x \cdot s)=\phi(x) \bullet f(s).
\end{eqnarray}
Furthermore, if $\phi$ and $f$ are bijections then we have a \emph{singquandle shadow isomorphism}.
\end{definition}

This definition implies that 
$(Im(f),Im(\phi),\triangleright,R'_1,R'_2,\bullet)$ is a \emph{subshadow} of $(S',Y,\triangleright,R'_1,R'_2,\bullet)$.\\


Let $L$ be an oriented singular link and $D$ be its diagram.  The set of arcs of $D$ will be denoted by $\mathcal{A}(D)$ and the set of connected regions of $\mathbb{R}^2\setminus D$ by $\mathcal{R}(D)$.  Definition~\ref{singHom} of a singquandle homomorphism allows to define the following notion of colorings by singquandles. 

\begin{definition}\cite{CCE}
Let $(S,*, R_1, R_2)$ be an oriented singquandle. An \emph{$S$-coloring} of $D$ is a map $f: \mathcal{A}(D) \rightarrow S$ such that at a crossing with  $u_1, u_2, o_1 \in \mathcal{A}(D)$ and at a singular crossing with $a_1,a_2,a_3,a_4 \in \mathcal{A}(D)$ the following conditions are satisfied,
\begin{eqnarray}
f(u_2) = f(u_1) * f(o_1), \\
f(a_3) = R_1(f(a_1),f(a_2)), \\
f(a_4) = R_2(f(a_1),f(a_2)).
\end{eqnarray}  
The conditions above are illustrated in Figure~\ref{shadowrule}.
\end{definition}

\begin{definition}\cite{CCE}
Let $(S,X,*, R_1, R_2, \cdot)$ be a shadow singquandle. An \emph{$(S,X)$-coloring} of $D$ is a map $f \times \phi: \mathcal{A}(D) \times \mathcal{R}(D) \rightarrow S \times X$ satisfying the following conditions,
\begin{itemize}
    \item $f$ is an $S$-coloring of $D$.
    \item $\phi(\mathcal{R}(D)) \subset X$.
    \item For $a \in \mathcal{A}(D)$ and $x_1, x_2 \in \mathcal{R}(D)$ the following
    \begin{equation}  
    \phi(x_1) \cdot f(a) = \phi(x_2).
    \end{equation}
\end{itemize}
The condition above is illustrated in Figure~\ref{shadowrule}.
\end{definition} 

When there is no confusion we will refer to an $(S,X)$-coloring by  a \emph{shadow coloring} of $D$.

The region coloring will be denoted by a box around the shadow element while the arc coloring by an element of a singquandle without a box.
Notice that the conditions required for the set $X$ to be an $S$-set for some oriented sinquandle, are the conditions needed to guarantee that shadow colorings are well defined at crossings, see Figure~\ref{shadowX}.

 The following results were introduced and proved in \cite{CCE}.

\begin{proposition}\label{countingshadow}
Let $L$ be a singular link diagram and $(S,X,*,R_1,R_2,\cdot)$ be a singquandle shadow. Then for each singquandle coloring of $L$ by $S$ and each element of $X$ there is exactly one shadow coloring of $L$.
\end{proposition}

For a singular link diagram $L$  and a singquandle shadow $(S,X)$. The \emph{shadow counting invariant}, $\#\textup{Col}_{(S,X)}(L)$, was defined in \cite{CCE} to be the number of shadow colorings of $L$ by $(S,X)$. From Proposition~\ref{countingshadow} and the definition of the shadow counting invariant one obtains the following result.

\begin{corollary}\label{singandshadow}
The shadow counting invariant of a singular link $L$ by the $S$-shadow $(S,X)$ is given by 
\[\# \textup{Col}_{(S,X)}(L) = \vert X \vert \, \# \textup{Col}_S(L) ,\]
where $\#\textup{Col}_S(L)$ is the singquandle counting invariant.
\end{corollary}
 This means that the shadow counting invariant is determined by the singquandle counting invariant. Therefore, the following polynomial was defined for a sinquandle shadow to extract more information from the singquandle shadow coloring set. 
\begin{definition}\cite{CCE}
The \emph{shadow singquandle polynomial}, denoted by $\textup{sp}(S,X)$, of the shadow singquandle $(S,X,*,R_1,R_2,\cdot)$ is the sum 
\[ \textup{sp}(S,X) = \sum_{x \in X} t^{r(x)},\]
where $r(x) = \vert \lbrace s \in S \, ; \, x \cdot s = x \rbrace\vert$. Furthermore, If $(S^{\prime},Y)$ is a subshadow of $(S,X)$, then the \emph{subshadow singquandle polynomial} of $(S^{\prime},Y)$ is 
\[ \textup{Subsp} ((S',Y) \subset (S,X)) = \sum_{x \in Y} t^{r(x)}, \]
where $r(x) =\vert \lbrace s' \in S' \, ; \, x \cdot s' = x \rbrace \vert$.
\end{definition}
The authors of \cite{CCE} were able to prove that the shadow singquandle polynomial is an invariant of singquandle shadows.

\begin{proposition} 
Let $(S,X)$ and $(S',Y)$ be two singquandle shadows.
If $(S,X)$ and $(S',Y)$  are isomorphic, then they have equal shadow polynomials, $\textup{sp}(S,X) = \textup{sp}(S',Y)$.
\end{proposition}


 The shadow polynomial can be used as a tool distinguish and classify singquandle shadows. The following example is taken from \cite{CCE}.
\begin{example}
Let $(S,*,R_1,R_2)$ be a singquandle with $S= \mathbb{Z}_8$, $x*y= 5x-4y = x \bar{*} y$, $R_1(x,y) = 3 x + 4 y$, and $R_2(x,y)= 4x+3y$. Let $X = \mathbb{Z}_4$. Consider the following singquandle shadow $(X, \cdot)$ with $\cdot s : X \rightarrow X$ for each $ s \in S$ defined by $x \cdot s = x + 2s +s^2$. The singquandle shadow $(S,X, *, R_1,R_2,\cdot)$ has the follow the following shadow polynomial 
\[  \textup{sp}(S,X) = 4t^4. \]
On the other hand, the singquandle shadow $(X, \bullet)$ with $\bullet s : X \rightarrow X$ for each $s \in S$ defined by $x \bullet s = 3x+2s +2s^2$. The singquandle shadow $(S,X, *, R_1,R_2,\bullet)$ has the following shadow polynomial 
\[ \textup{sp}(S,W) = 2+2t^8. \]
\end{example}

The subshadow singquandle polynomial and the shadow coloring invariant were combined in order to produce an invariant of singular links. This polynomial invariant detected different information than the shadow counting invariant and that the previous defined singquandle polynomial invariant.

\begin{definition}\cite{CCE}
Let $f \times \phi$ be a shadow coloring of an oriented singular link diagram $D$. The closure of the set of shadow colors under the action of the image subsingquandle $\textup{Im}(f) \subset S$ of $f \times \phi$ is a subshadow called the \emph{shadow image} of $f \times \phi$, which we denote by $\textup{om}(f \times \phi).$ 
\end{definition}


\begin{definition}
Let $(S,X)$ be an $S$-shadow and let $L$ be an oriented singular link with diagram $D$. The \emph{singquandle shadow polynomial invariant} of $L$ with respect $(S,X)$ is 
\[ SP(L) = \sum_{f \times \phi \in \textup{shadow coloring}} u^{\textup{Subsp}  (\textup{om}(f \times \phi) \subset (S,X))}.\]
\end{definition}

The following example from \cite{CCE} was used to test the strength of the singquandle shadow polynomial invariant of $L$. Precisely, the following pair of singular knots have the same number of singquandle colorings and the same singquandle polynomial invariant, but are distinguished by their singquandle shadow polynomial invariant.

\begin{example}
Let $(S,X,*,R_1,R_2,\cdot)$ be the shadow singquandle with $S=\mathbb{Z}_8$, $X=\mathbb{Z}_6$ and operations $x*y = 3x-2y = x \bar{*} y$, $R_1(x,y) = 7x+6y$, $R_2(x,y) = 2x+3y$, and shadow operation $\cdot s : X \rightarrow X$ for all $s \in S$ defined by $x \cdot s = x + 3 s.$

The pair of singular knots in Figure~\ref{Example1} have the same singquandle counting invariant $\#\textup{Col}_S(4_1^k)  = 16 =  \#\textup{Col}_S(5_4^k)$. Therefore, by Theorem~\ref{singandshadow} we obtain that the two singular knots have the same shadow counting invariant $\#\textup{Col}_{(S,X)}(4_1^k)  = 96 = \#\textup{Col}_{(S,X)}(5_4^k)$. Furthermore, the two singular knots have the same singquandle polynomial $\phi_{Ssqp}(4_1^k)=4 u^{s_1^2 t_1^2 s_2^2 t_2^2 s_3 t_3}+4 u^{2 s_1^2 t_1^2 s_2^2 t_2^2 s_3
   t_3}+8 u^{4 s_1^2 t_1^2 s_2^2 t_2^2 s_3 t_3} =\phi_{Ssqp}(5_4^k) $. However, the singquandle shadow polynomial invariant distinguishes the two singular knots:

\begin{figure}[h]
    \centering

\begin{tikzpicture}[use Hobby shortcut,scale=1,add arrow/.style={postaction={decorate}, decoration={
  markings,
  mark=at position .5 with {\arrow[scale=1.5,>=stealth]{<}}}}]
\begin{knot}[
  consider self intersections=true,
clip width=3,
  flip crossing/.list={2,4,6}
]
\strand ([closed]90:2) foreach \k in {1,...,5} { .. (90-360/5+\k*720/5:1.5) .. (90+\k*720/5:2) } (90:2)[add arrow];
\end{knot}

\node[circle,draw=black, fill=black, inner sep=0pt,minimum size=6pt] (a) at (1,1.35) {};

\node[above] at (0,1) {\tiny $s_1$};
\node[above] at (0,2) {\tiny $s_2$};
\node[left] at (1.5,.5) {\tiny $s_3$};
\node[right] at (2,.5) {\tiny $s_4$};
\node[left] at (-1.5,-1.5) {\tiny $s_5$};
\node[left] at (-2,.5) {\tiny $s_6$};

\draw (0,0) node [draw] {\tiny $x$};

\node[below] at (0,-2.5){$SP(4_1^k)= 24 u^{t^2}+24 u^t+48 u^2$};

\end{tikzpicture}
\hspace{1in}
\begin{tikzpicture}[use Hobby shortcut,scale=1,add arrow/.style={postaction={decorate}, decoration={
  markings,
  mark=at position .5 with {\arrow[scale=1.5,>=stealth]{>}}}}]
\begin{knot}[
  consider self intersections=true,
  ignore endpoint intersections=false,
 flip crossing/.list={7,4,9},
  clip width=4,
  only when rendering/.style={
  }
  ]
\strand ([closed]0,2)..(-.5,1.75)..(0,0)..(.5,-.2)..(2,-2)..(-.5,-1)..(-2,1.5)..(0,1.5)..(2,1.5)..(.5,-1)..(-2,-2)..(-.5,-.2)..(0,0)..(.5,1.75)..(0,2)[add arrow];

\end{knot}

\node[circle,draw=black, fill=black, inner sep=0pt,minimum size=6pt] (a) at (0,0) {};
\node[right] at (-1,0) {\tiny $s_1$};
\node[right] at (-1,.5) {\tiny $s_2$};
\node[right] at (.5,0) {\tiny $s_3$};
\node[right] at (.5,.5) {\tiny $s_4$};
\node[left] at (-2,.5) {\tiny $s_5$};
\node[right] at (2,.5) {\tiny $s_6$};
\node[left] at (-2,-1) {\tiny $s_7$};

\draw (0,-.5) node [draw] {\tiny $x$};

\node[below] at (0,-2.5){$SP(5_4^k)=48 u^{t^4}+24 u^{t^2}+24 u^t $};

\end{tikzpicture}    
    \caption{Singular knots $4_1^k$ and $5_4^k$ and corresponding $SP$ invariant.}
    \label{Example1}
\end{figure}
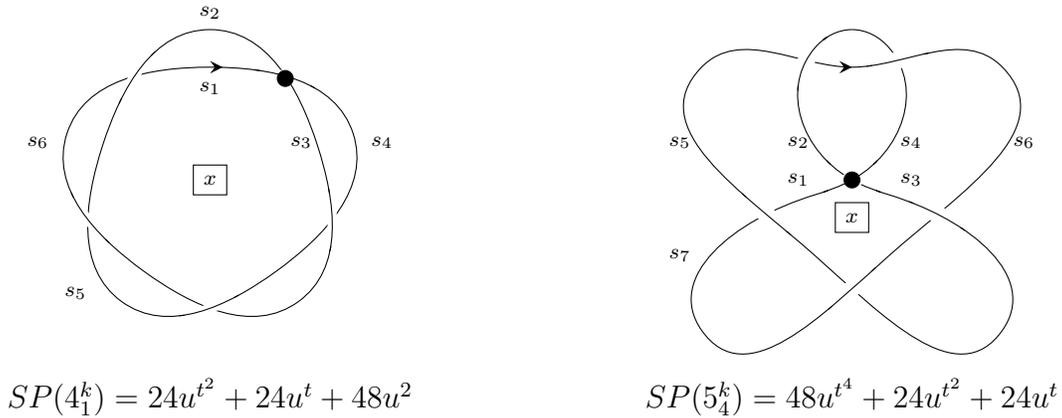
\end{example}

\section{Psyquandles and Singular Link Invariants}\label{Psy}

In \cite{NOS}, the authors introduced the notion of \emph{psyquandle} as an algebraic generalization of the notion of biquandle that entails the topological moves characterizing both singular links and pseudolinks.  The authors showed that the corresponding notion of psyquandle colorings of an oriented singular link, or pseudolink, determines a counting invariant. Furthermore, in \cite{CN2}, the authors gave a cocycle enhancement theory to the case of psyquandles. Precisely, they defined enhancements of the psyquandle counting invariant via pairs of a biquandle 2-cocycle and a new function satisfying some conditions. As an application they defined new single-variable and two-variable polynomial invariants of oriented pseudoknots and singular knots and links. In this section we will provide an overview of results obtained for signular knots from psyquandles. The following definition was initially introduced in \cite{NOS}, but was reformulated in \cite{CN2}. We will present the reformulated version of the definition.

\begin{definition}\cite{CN2}
Let $X$ be a set. A \textit{psyquandle structure} on $X$ is set of four
binary operations $\utr,\otr,\ud,\od:X\times X\to X$ satisfying the conditions

(I) All four operations are right-invertible,i.e. there exist
binary operations $\utr^{-1},\otr^{-1},\ud^{-1},\od^{-1}:X\times X\to X$ such that
\[\begin{array}{rcccl}
(x\utr y)\utr^{-1}y & = & (x\utr^{-1}y)\utr y & = & x\\ 
(x\otr y)\otr^{-1}y & = & (x\otr^{-1}y)\otr y & = & x\\ 
(x\ud y)\ud^{-1}y & = & (x\ud^{-1}y)\ud y & = & x\\ 
(x\od y)\od^{-1}y & = & (x\od^{-1}y)\od y & = & x,
\end{array}\]

(II) For all $x\in X$, $x\utr y=x\otr y$, 

(III) For all $x,y\in X$, the maps $S,S':X\times X\to X\times X$ 
defined by
\[S(x,y)=(y\otr x,x\utr y) \quad\mathrm{and}\quad S'(x,y)=(y\od x,x\ud y)\] 
are invertible,

(IV) For all $x,y,z\in X$,
\[\begin{array}{rcl}
(x\utr y)\utr (z\utr y) & = & (x\utr z)\utr (y\otr z) \\
(x\utr y)\otr (z\utr y) & = & (x\otr z)\utr (y\otr z) \\
(x\otr y)\otr (z\otr y) & = & (x\otr z)\otr (y\utr z) \\
\end{array}\]

(V) For all $x,y\in X$ we have
\[\begin{array}{rcl}
x\ud((y\otr x)\od^{-1} x) & = & [(x\utr y)\od^{-1} y]\otr[(y\otr x)\ud^{-1} x]\\
y\ud((x\utr y)\od^{-1} y) & = & [(y\otr x)\od^{-1} x]\utr[(x\utr y)\od^{-1} y],
\end{array}\]
and

(VI) For all $x,y,z\in X$ we have
\[\begin{array}{rcl}
(x\otr y)\otr (z\od y) & = & (x\otr z)\otr (y\ud z) \\
(x\utr y)\utr (z\od y) & = & (x\utr z)\utr (y\ud z) \\
(x\otr y)\od (z\otr y) & = & (x\od z)\otr (y\utr z) \\
(x\utr y)\ud (z\utr y) & = & (x\ud z)\utr (y\otr z) \\
(x\otr y)\ud (z\otr y) & = & (x\ud z)\otr (y\utr z) \\
(x\utr y)\od (z\utr y) & = & (x\od z)\utr (y\otr z).
\end{array}\]

A psyquandle which also satisfies $x\ud x=x\od x$ for all
$x\in X$ is said to be \textit{pI-adequate}. 
\end{definition}

\begin{figure}[h]
\begin{tikzpicture}[use Hobby shortcut]
\begin{knot}[
consider self intersections=true,
  ignore endpoint intersections=false,
  flip crossing/.list={3,4,5,6,11,13},
  clip width = 4
]
\strand[decoration={markings,mark=at position 1 with
    {\arrow[scale=3,>=stealth]{>}}},postaction={decorate}] (1,1) ..(-1,-1); 
\strand[decoration={markings,mark=at position 1 with
    {\arrow[scale=3,>=stealth]{>}}},postaction={decorate}] (-1,1) ..(1,-1);
\strand[decoration={markings,mark=at position 1 with
    {\arrow[scale=3,>=stealth]{>}}},postaction={decorate}] (2,1) ..(4,-1);
\strand[decoration={markings,mark=at position 1 with
    {\arrow[scale=3,>=stealth]{>}}},postaction={decorate}] (4,1) ..(2,-1);  
\strand[decoration={markings,mark=at position 1 with
    {\arrow[scale=3,>=stealth]{>}}},postaction={decorate}] (5,1) ..(7,-1);
\strand[decoration={markings,mark=at position 1 with
    {\arrow[scale=3,>=stealth]{>}}},postaction={decorate}] (7,1) ..(5,-1);      
\end{knot}
\node[left] at (-1,1) {\tiny $x$};
\node[left] at (-1,-1) {\tiny $y$};
\node[left] at (.8,1) {\tiny $y \otr x$};
\node[left] at (.8,-1) {\tiny $ x \utr y $};

\node[left] at (2,1) {\tiny $x$};
\node[left] at (2,-1) {\tiny $y$};
\node[left] at (3.8,1) {\tiny $x \utr y$};
\node[left] at (3.8,-1) {\tiny $y \otr x$};

\node[circle,draw=black, fill=black, inner sep=0pt,minimum size=8pt] (a) at (6,0) {};
\node[left] at (5,1) {\tiny $x$};
\node[left] at (5,-1) {\tiny $y$};
\node[left] at (6.8,1) {\tiny $y \od x$};
\node[left] at (6.8,-1) {\tiny $x \ud y$};
\end{tikzpicture}
\vspace{.2in}
		\caption{Semiarc coloring rule at classical and singular crossings.}
		\label{semicol}
\end{figure}
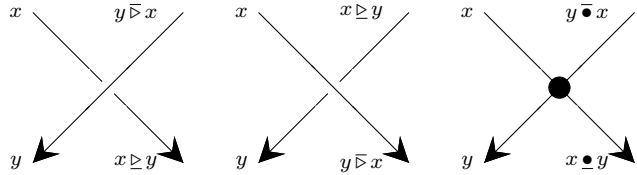

A pI-adaquate psyquandle is an invariant of pseudolinks. For more details and results regarding pI-adaquate psyquandles the reader is referred to \cite{NOS, CN2}. In this paper we will focus on the results for singular links obtained from psyquandles. The axioms of a psyquandle structure are motivated by the generalized Reidemeister moves for singular knots and  the semi-arc coloring rules in Figure~\ref{semicol}. This motivation is different from the singquandle axioms which are motivated by the arc coloring rules in Figure~\ref{ColSing}.

The following examples and additional examples can be found in \cites{NOS,CN2}.

\begin{example}
Every biquandle is a psyquandle by setting $x\od y=x\otr y$ and
$x\ud y=x\utr y$. 
\end{example}
\begin{example}\label{psy}
The psyquandle structure on the set $X=\{1,2,3,4,5,6\}$, can be specified by explicitly listing the operation tables of $\utr, \otr, \ud, \od$. In practice it is convenient to put these together 
into an $6\times 24$ bock matrix, so the psyquandle structure on $X=\{1,2,3,4,5,6\}$
specified by
\[\left[\begin{array}{cccccc|cccccc|cccccc|cccccc}
2& 4 & 4 & 6 & 6 & 2 &  2 & 6 & 2 & 6 & 2 & 6 &  2& 4& 2& 6& 2& 2 &  2& 6& 4& 6& 6& 6\\ 
3& 5 & 5 & 1 & 1 & 3 &  1 & 5 & 1 & 5 & 1 & 5 &  3& 5& 5& 5& 1& 5 &  1& 5& 1& 1& 1& 3\\
4& 6 & 6 & 2 & 2 & 4 &  6 & 4 & 6 & 4 & 6 & 4 &  6& 6& 6& 2& 6& 4 &  4& 4& 6& 4& 2& 4\\
5& 1 & 1 & 3 & 3 & 5 &  5 & 3 & 5 & 3 & 5 & 3 &  5& 3& 1& 3& 3& 3 &  5& 1& 5& 3& 5& 5\\
6& 2 & 2 & 4 & 4 & 6 &  4 & 2 & 4 & 2 & 4 & 2 &  4& 2& 4& 4& 4& 6 &  6& 2& 2& 2& 4& 2\\
1& 3 & 3 & 5 & 5 & 1 &  3 & 1 & 3 & 1 & 3 & 1 &  1& 1& 3& 1& 5& 1 &  3& 3& 3& 5& 3& 1
\end{array}\right].\]
\end{example}

Furthermore, the psyquandle structure was used to define an invariant of singular knots called \emph{psyquandle counting invariant}, denoted by $\Phi_X^\mathbb{Z}(L)$, where $X$ is a finite psyquandle, $L$ is an oriented singular knot or link in \cite{NOS}. The psyquandle counting invariant is the cardinality of $\mathcal{C}(X,L)$, the set of colorings of the singular link $L$ by a finite psyquandle $X$. The authors of \cite{NOS,CN2} refer to an element of $\mathcal{C}(L,X)$ as an $X$-coloring of $L$. The authors of \cite{CN2} defined an enhancement of the psyquandle counting invariant by generalizing the biquandle cocycle invariant to the case of psyquandles. 

\begin{definition}\cite{CN2}\label{CocyclesDef}
Let $X$ be a psyquandle and $R$ a commutative ring with identity. A 
\textit{Boltzmann weight} for $X$ is a pair of maps $\phi,\psi:X\times X\to R$
satisfying

(I) For all $x\in X$, $\phi(x,x)=0$

(II) For all $x,y\in X$,
\[\phi(x,y)+\psi(y,(x\utr y)\od^{-1}y)=\phi((y\otr x)\od^{-1} x, (x\utr y)\od^{-1}y)+\psi(x,(y\otr x)\od^{-1} x).\]

(III) For all $x,y, z \in X$,
\begin{eqnarray*}
\phi(x,y)+\phi(y,z)+\phi(x\utr y,z\otr y) & = & \phi(x\utr z,y\utr z)+\phi(x,z) + \phi(y\otr x,z\otr x)\\
\psi(x,y)+\phi(y,z)+\phi(x \ud y, z\otr y) &=& \psi(x \utr z, y \utr z)+\phi(x,z) + \phi(y \od x, z \otr x) \\
\psi (z,y) - \phi(x,y) - \phi(x \utr y, z \ud y) &=& \psi(z \otr x, y \otr x) -\phi(x,z) - \phi(x \utr z, y \od z).
\end{eqnarray*}

We say that $\phi$ and $\psi$ are \textit{strongly compatible} if
we also have

(IV) For all $x,y,z\in X$,
\[\psi(x,y)=\psi(x\utr z, y\utr z)\quad\mathrm{and}\quad 
\psi(z,y)=\psi(z\otr x, y\otr x).\]

\end{definition}

\begin{figure}[h]
\begin{tikzpicture}[use Hobby shortcut]
\begin{knot}[
consider self intersections=true,
  ignore endpoint intersections=false,
  flip crossing/.list={3,4,5,6,11,13},
  clip width = 4
]
\strand[decoration={markings,mark=at position 1 with
    {\arrow[scale=3,>=stealth]{>}}},postaction={decorate}] (1,1) ..(-1,-1); 
\strand[decoration={markings,mark=at position 1 with
    {\arrow[scale=3,>=stealth]{>}}},postaction={decorate}] (-1,1) ..(1,-1);
\strand[decoration={markings,mark=at position 1 with
    {\arrow[scale=3,>=stealth]{>}}},postaction={decorate}] (2,1) ..(4,-1);
\strand[decoration={markings,mark=at position 1 with
    {\arrow[scale=3,>=stealth]{>}}},postaction={decorate}] (4,1) ..(2,-1);  
\strand[decoration={markings,mark=at position 1 with
    {\arrow[scale=3,>=stealth]{>}}},postaction={decorate}] (5,1) ..(7,-1);
\strand[decoration={markings,mark=at position 1 with
    {\arrow[scale=3,>=stealth]{>}}},postaction={decorate}] (7,1) ..(5,-1);      
\end{knot}
\node[above] at (-1,1) {\tiny $x$};
\node[above] at (1,1) {\tiny $y$};
\node[below] at (0,-1) {\tiny $+\phi(x,y)$};

\node[below] at (2,-1) {\tiny $x$};
\node[below] at (4,-1) {\tiny $y$};
\node[below] at (3,-1) {\tiny $-\phi(x,y)$};

\node[circle,draw=black, fill=black, inner sep=0pt,minimum size=8pt] (a) at (6,0) {};
\node[above] at (5,1) {\tiny $x$};
\node[above] at (7,1) {\tiny $y$};
\node[below] at (6,-1) {\tiny $+\psi(x,y)$};
\end{tikzpicture}
\vspace{.2in}
		\caption{Boltzmann weights at classical and singular crossings.}
		\label{Bweights}
\end{figure}
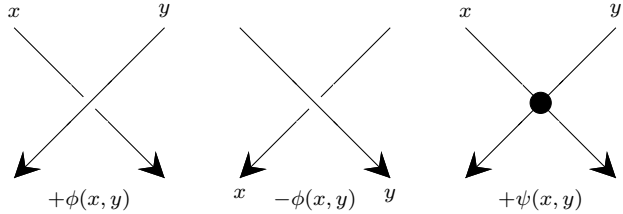

The Boltzmann weight axioms are motivated by generalized Reidemeister moves and using the contribution rule in Figure~\ref{Bweights}.

\begin{definition}\cite{CN2}
Let $X$ be a psyquandle, $R$
a commutative ring with identity and $(\phi,\psi)$ a Boltzmann weight. Let $L$ be an oriented 
singular knot or link.

(I) For each $X$-coloring $L_c$ of $L$, we define the \textit{Boltzmann weight} of $L_c$,
denoted $BW(L_c)$, to be the sum of contributions over all crossings in $L_c$.

(II)
We define the \textit{single-variable Boltzmann-enhanced psyquandle polynomial}
to be 
\[\Phi_X^{\phi,\psi}(L)\sum_{L_c\in\mathcal{C}(L,X)} w^{BW(L_c)}.\]

(III) If $\phi$ and $\psi$ are strongly compatible, we define the 
\textit{partial Boltzmann weights} $BW_\phi(L_c)$ and $BW_\psi(L_c)$ to be 
the sums of $\phi$ contributions and $\psi$ contributions respectively; then
we define the \textit{two-variable Boltzmann-enhanced psyquandle polynomial}
to be
\[\Phi_X^{\phi,\psi}(L)\sum_{L_c\in\mathcal{C}(L,X)} u^{BW_{\phi}(L_c)}v^{BW_{\psi}(L_c)}.\] 

\end{definition}

Furthermore, both the single-variable Boltzmann-enhanced psyquandle polynomial and the two-variable Boltzmann-enhanced psyquandle polynomial were shown to be an invariant of singular links. The following example and other examples can be found in \cite{CN2}. A link to the custom \texttt{Python} code used to compute the examples can be found in \cite{CN2}. 
\begin{example}\label{exlist}
In this example the single-variable Boltzmann-enchanced psyquandle polynomial for the \emph{2-bouquet graphs of type $L$} (with choice of orientation) in \cite{Oyamaguchi} are collected. The psyquandle structure in Example~\ref{psy} was used along 
with Boltzmann weight function $\phi: X \times X \rightarrow \mathbb{Z}_2$ given by the matrix
\[\left[\begin{array}{cccccc}
0& 1& 0& 1& 0& 1\\
0& 0& 0& 0& 0& 0\\
0& 1& 0& 1& 0& 1\\
0& 0& 0& 0& 0& 0\\
0& 1& 0& 1& 0& 1\\
0& 0& 0& 0& 0& 0
\end{array}\right]\]
and weight function $\psi: X \times X \rightarrow \mathbb{Z}_2$ given by the matrix
\[\left[\begin{array}{cccccc}
1& 0& 1& 0& 1& 0\\
1& 1& 1& 1& 1& 1\\
1& 0& 1& 0& 1& 0\\
1& 1& 1& 1& 1& 1\\
1& 0& 1& 0& 1& 0\\
1& 1& 1& 1& 1& 1
\end{array}\right].\]
The results are collected in the table

\begin{center}
\begin{tabular}{ l| c| c }
 2-bouquet graph of type $L$ &$\Phi_X^\mathbb{Z}(L)$ & $\Phi_X^{\phi,\psi}(L)$  \\
\hline
$5_2^l,6_1^l$                   &$12$       & $12w$ 	 \\
 \hdashline
$3_1^l,4_1^l,5_3^l,6_2^l,6_6^l$ &           & $6w+6$	\\
\hline
$6_3^l,6_8^l,6_9^l,6_{10}^l,6_{11}^l$   &   $24$ 	& $24w$ 	\\	
 \hdashline
$5_1^l,6_5^l,6_7^l$ & 	    & $6w+18$ 	   	 \\
 \hdashline
$1_1^l$ &  		& $18w+6$ 	  	\\
\hline
$6_4^l, 6_{12}^l$ & $36$    & $18w+18$.	
\end{tabular}
\end{center}
\end{example}

{\bf Acknowledgements}
 Mohamed Elhamdadi was partially supported by Simons Foundation collaboration grant 712462.

\end{document}